# Optimal rates for plug-in estimators of density level sets

PHILIPPE RIGOLLET[1] and RÉGIS VERT[2]

[1]*Department of Operations Research and Financial Engineering, Princeton University, Princeton, NJ 08544, USA. E-mail: rigollet@princeton.edu*
[2]*Masagroup, 24 Bd de l'Hôpital, 75005 Paris, France. E-mail: regis.vert@masagroup.net*

In the context of density level set estimation, we study the convergence of general plug-in methods under two main assumptions on the density for a given level $\lambda$. More precisely, it is assumed that the density (i) is smooth in a neighborhood of $\lambda$ and (ii) has $\gamma$-exponent at level $\lambda$. Condition (i) ensures that the density can be estimated at a standard nonparametric rate and condition (ii) is similar to Tsybakov's margin assumption which is stated for the classification framework. Under these assumptions, we derive optimal rates of convergence for plug-in estimators. Explicit convergence rates are given for plug-in estimators based on kernel density estimators when the underlying measure is the Lebesgue measure. Lower bounds proving optimality of the rates in a minimax sense when the density is Hölder smooth are also provided.

## 1. Introduction

Let $Q$ be a positive $\sigma$-finite measure on $\mathcal{X} \subseteq \mathbb{R}^d$. Consider i.i.d. random vectors $(X_1, \ldots, X_n)$ with distribution $P$, having an unknown probability density $p$ with respect to the measure $Q$. For a fixed $\lambda > 0$, we are interested in the estimation of the $\lambda$-*level set* of the density $p$ defined by

$$\Gamma_p(\lambda) \triangleq \{x \in \mathcal{X} : p(x) > \lambda\}. \tag{1.1}$$

Throughout the paper, we fix $\lambda > 0$ and when no confusion is possible, we use the notation $\Gamma(\lambda)$, or simply $\Gamma$, instead of $\Gamma_p(\lambda)$. When $Q$ is the Lebesgue measure on $\mathbb{R}^d$, density level sets typically correspond to minimum volume sets of given $P$-probability mass, as shown in Polonik (1997).

**Remark 1.1.** A somewhat preponderant definition of a density level set is

$$\overline{\Gamma}(\lambda) \triangleq \{x \in \mathcal{X} : p(x) \geq \lambda\}, \tag{1.2}$$

that is, the union of $\Gamma(\lambda)$ and the set $\{x \in \mathcal{X} : p(x) = \lambda\}$. Since, in this paper, the density is allowed to have flat parts at level $\lambda$, the sets $\Gamma(\lambda)$ and $\overline{\Gamma}(\lambda)$ can differ by an arbitrarily large set. Density level sets defined by (1.1) or (1.2) can be estimated using plug-in estimators with positive or negative offset, respectively (see Section 2.2). However, definition (1.1) remains consistent with the definition of the support of the density when $\lambda = 0$. The results detailed hereafter pertain only to this definition, but are applicable to definition (1.2) after minor changes.

The following are two possible applications of density level set estimation.

**Anomaly detection.** The goal is to detect an abnormal observation from a sample (see, e.g., Steinwart *et al.* (2005) and references therein). One way to deal with that problem is to assume that abnormal observations do not belong to a group of concentrated observations. In this framework, observations are considered to be abnormal when they do not belong to $\Gamma(\lambda)$ for some fixed $\lambda \geq 0$. The special case $\lambda = 0$, which corresponds to support estimation, has been examined by Devroye and Wise (1980), for example. In the general case, $\lambda$ can be considered a tolerance level for anomalies: the smaller $\lambda$, the fewer observations are considered to be abnormal.

**Unsupervised or semi-supervised classification.** These two problems amount to the identification of areas where the observations are concentrated with possible use of some available labels for the semi-supervised case. For instance, it can be assumed that the connected components of $\Gamma(\lambda)$, for a fixed $\lambda$, are clusters of homogeneous observations, as described in Hartigan (1975). Note that this definition has been refined, for example in Stuetzle (2003), and has been studied with plug-in estimators in Rigollet (2007).

*Remark 1.2.* In both applications, the choice of $\lambda$ is critical and must be addressed carefully. However, this problem is beyond the scope of this paper.

There are essentially two approaches toward estimating density level sets from the sample $(X_1, \ldots, X_n)$. The most straightforward is to resort to *plug-in* methods where the density $p$ in the expression for $\Gamma_p(\lambda)$ is replaced by its estimate computed from the sample. Another way to estimate density level sets is to resort to *direct* methods which are based on empirical excess-mass maximization. The *excess-mass H* is a functional that measures the quality of an estimator $\hat{G}$ and is defined as follows Hartigan (1987); Müller and Sawitzki (1987):

$$H(\hat{G}) = P(\hat{G}) - \lambda Q(\hat{G}).$$

Excess-mass measures how the $P$-probability mass concentrates in the region $\hat{G}$ and it is maximized by $\Gamma = \Gamma(\lambda)$. Hence, it acts as a risk functional in the density level set estimation (DLSE) framework and it is natural to measure the performance of an estimator $\hat{G}$ by its *excess-mass deficit* $H(\Gamma) - H(\hat{G}) \geq 0$. Further justifications for the well-foundedness of the excess-mass criterion can be found in Polonik (1995). Recently,

Gayraud and Rousseau (2005) proposed a Bayesian approach to DLSE, together with interesting comparative simulations.

While local versions of direct methods have been analyzed deeply and proven to be optimal in a minimax sense over a certain family of well-behaved distributions (see Tsybakov (1997)), and although reasonable implementations have been recently proposed (see, e.g., Steinwart *et al.* (2005)), they are still not very easy to use for practical purposes, compared to plug-in methods. Indeed, in practice, rather than specifying a value for $\lambda$, the user can specify a value for $\alpha$, the $P$-probability mass of the level set. In this case, the value of $\lambda$ is implied by that of $\alpha$ and efficient direct methods can be derived (Scott and Nowak (2006)). However, in general, using direct methods, one must run an optimization procedure several times, for different density level values, then choose a posteriori the most suitable level according to the desired rejection rate. Plug-in methods do not involve such a complex process: the density estimation step is only performed once and the construction of a density level set estimate simply amounts to thresholding the density estimate at the desired level.

On the other hand, in the related context of binary classification, where more theoretical advances have been developed, the different analyses proposed thus far have mainly supported a belief in the superiority of direct methods. Yang (1999) shows that, under general assumptions, plug-in estimators cannot achieve a classification error risk convergence rate faster than $O(1/\sqrt{n})$ – where $n$ is the size of the data sample – and suffer from the curse of dimensionality. In contrast to that, under slightly different assumptions, direct methods achieve this rate $O(1/\sqrt{n})$, whatever the dimensionality (see, e.g., Vapnik (1998); Devroye *et al.* (1996); Tsybakov (2004)), and can even reach faster convergence rates – up to $O(1/n)$ – under *Tsybakov's margin assumption* (see Mammen and Tsybakov (1999); Tsybakov (2004); Tsybakov and van de Geer (2005); Tarigan and van de Geer (2006)). This contributed to the arousing of some pessimism concerning plug-in methods. Nevertheless, such a comparison between plug-in methods and direct methods is far from legitimate, since the aforementioned analyses of both plug-in methods and direct ones have been carried out under different sets of assumptions (those sets are not disjoint, but none of them is included in the other).

Recently, in the standard classification framework, Audibert and Tsybakov (2007) have combined a new type of assumption dealing with the smoothness of the *regression function* and the well-known margin assumption. Under these assumptions, they derive fast convergence rates – even faster than $O(1/n)$ in some situations – for plug-in classification rules based on local polynomial estimators. This new result reveals that plug-in methods should not be considered inferior to direct methods and, more importantly, that this new type of assumption on the regression function is a critical point in the general analysis of classification procedures.

In this paper, we extend such positive results to the DLSE framework: we revisit the analysis of plug-in density level set estimators and show that they can be also very efficient under smoothness assumptions on the underlying density function $p$. Unlike the global smoothness assumption used in Audibert and Tsybakov (2007), the local smoothness assumption introduced here emphasizes the predominant role of the smoothness close to the level $\lambda$, as opposed to the smoothness for values of $p$ far from the level under consideration. Related papers are Baíllo *et al.* (2001) and Baíllo (2003), which investigate plug-in estimators based on a certain type of kernel density estimate. Baíllo *et al.* (2001) also study the convergence for the symmetric difference under other assumptions and Baíllo (2003) derives almost sure rates of convergence for a quantity different from the ones studied here. It is interesting to observe that she introduces a condition similar to the $\gamma$-exponent used here.

The particular case $\lambda = 0$ corresponds to estimation of the support of density $p$ and is often applied to anomaly detection. Following the pioneering paper of Devroye and Wise (1980), this problem has received more attention than the general case $\lambda \geq 0$ and has been treated using plug-in methods, for example by Cuevas and Fraiman (1997). Unlike the previously cited papers, we derive rates of convergence and prove that these rates are optimal in a minimax sense. However, we do not treat the case $\lambda = 0$ for the which the rates are typically different than for $\lambda > 0$, as pointed out by Tsybakov (1997), for example. The techniques employed in the present analysis can be refined to encompass the case $\lambda = 0$ and the results will be published separately.

A general plug-in approach has been studied previously by Molchanov (1998), where a result on the asymptotic distribution of the Hausdorff distance is given. In a recent paper, Cuevas *et at.* (2006) study general plug-in estimators of the level sets. Under very general assumptions, they derive consistency with respect to the Hausdorff metric and the measure of the symmetric difference. However, this very general framework does not allow them to derive rates of convergence.

This paper is organized as follows. Section 2 introduces the notation and definitions. Section 3 presents the main result, that is, a new bound on the error of plug-in estimators based on general density estimators that satisfy a certain exponential inequality. In Section 4, we then apply this result to the particular case of kernel density estimators, under the assumption that the underlying density belongs to some locally Hölder smooth class of densities. Finally, minimax lower bounds are given in Section 5, as a way to assess the optimality of the upper bounds involved in the main result.

## 2. Notation and setup

For any vector $x \in \mathbb{R}^d$, denote by $x^{(j)}$ its $j$th coordinate, $j = 1, \ldots, d$. Denote by $\|\cdot\|$ the Euclidean norm in $\mathbb{R}^d$ and by $\mathcal{B}(x, r)$ the closed Euclidean ball in $\mathcal{X}$ centered at $x \in \mathcal{X}$ and of radius $r > 0$.

The probability and expectation with respect to the joint distribution of $(X_1, \ldots, X_n)$ are denoted by $\mathbb{P}$ and $\mathbb{E}$, respectively. For any function $f : \mathbb{R}^d \to \mathbb{R}$, we denote by $\|f\|_\infty = \sup_{x \in \mathbb{R}^d} |f(x)|$ the sup-norm of $f$ and by $\|f\| = (\int_{\mathbb{R}^d} f^2(x) \,\mathrm{d}x)^{1/2}$ its $L_2$-norm. Also, for any measurable function $f$ on $\mathcal{X}$ and any set $A \subset f(\mathcal{X})$, we write, for simplicity, $\{x \in \mathcal{X} : f(x) \in A\} = \{f \in A\}$. Throughout the paper, we denote by $C$ positive constants that can change from line to line and by $c_j$ positive constants that have to be identified. Finally, $A^c$ denotes the complement of the set $A$.

Our choice of measure of performance will affect the construction of our estimator, so we begin by discussing this topic.

### 2.1. Measures of performance

Recall that $Q$ is a positive $\sigma$-finite measure on $\mathcal{X}$ and define the measure $\tilde{Q}_\lambda$ that has density $|p(\cdot) - \lambda|$ with respect to $Q$. To assess the performance of a density level set estimator, we use the two pseudo-distances between two sets $G_1$ and $G_2 \subseteq \mathcal{X}$:

(i) the $Q$-measure of the symmetric difference between $G_1$ and $G_2$,

$$d_\Delta(G_1, G_2) = Q(G_1 \Delta G_2);$$

(ii) the $\tilde{Q}_\lambda$-measure of the symmetric difference between $G_1$ and $G_2$,

$$d_H(G_1, G_2) = \tilde{Q}_\lambda(G_1 \Delta G_2) = \int_{G_1 \Delta G_2} |p(x) - \lambda| \, dQ(x).$$

The quantity $d_\Delta(G_1, G_2)$ is a standard and natural way to measure the distance between two sets $G_1$ and $G_2$. Note that for any measurable set $G \subseteq \mathcal{X}$, the excess-mass $H(G)$ can be written

$$H(G) = \int_G (p(x) - \lambda) \, dQ(x).$$

Thus, we can rewrite

$$H(\Gamma) - H(\hat{G}) = \int_{\mathcal{X}} (\mathbb{1}_{\{p(\cdot) \geq \lambda\}}(x) - \mathbb{1}_{\hat{G}}(x))(p(x) - \lambda) \, dQ(x)$$

$$= \int_{\Gamma \Delta \hat{G}} |p(x) - \lambda| \, dQ(x) = d_H(\hat{G}, \Gamma).$$

This explains the notation $d_H$.

The following definition introduces a quantity which critically controls the complexity of the problem and therefore the attainable rates of convergence. In particular, it allows us to link $d_H$ to $d_\Delta$.

**Definition 2.1.** *For any $\lambda, \gamma \geq 0$, a function $f : \mathcal{X} \to \mathbb{R}$ is said to have $\gamma$-**exponent at level** $\lambda$ with respect to $Q$ if there exist constants $c_0 > 0$ and $\varepsilon_0 > 0$ such that, for all $0 < \varepsilon \leq \varepsilon_0$,*

$$Q\{x \in \mathcal{X} : 0 < |f(x) - \lambda| \leq \varepsilon\} \leq c_0 \varepsilon^\gamma.$$

The assumption under which the underlying density has $\gamma$-exponent at level $\lambda$ was first introduced by Polonik (1995). Its counterpart in the context of binary classification is commonly referred to as *margin assumption* (see Mammen and Tsybakov (1999); Tsybakov (2004)).

The exponent $\gamma$ controls the slope of the function around level $\lambda$. When $\gamma = 0$, the condition holds trivially and when $\gamma$ is positive, it constrains the rate at which the function approaches the level $\lambda$. A standard case corresponds to $\gamma = 1$, arising, for instance, in the case where the gradient of $f$ has a coordinate bounded away from 0 in a neighborhood of $\{f = \lambda\}$.

We now show that the pseudo-distances $d_\Delta$ and $d_H$ are linked when the density $p$ has $\gamma$-exponent at level $\lambda$. The following proposition is a direct consequence of Proposition A.1.

**Proposition 2.1.** *Fix $\lambda > 0$ and $\gamma \geq 0$. If the density $p$ has $\gamma$-exponent at level $\lambda$ w.r.t. $Q$, then, for any $L_Q > 0$, there exists $C > 0$ such that for any $G_1, G_2$ satisfying $Q(G_1 \Delta G_2) \leq L_Q$, we have*

$$d_\Delta(G_1, G_2) \leq Q(G_1 \Delta G_2 \cap \{p = \lambda\}) + C(d_H(G_1, G_2))^{\gamma/(1+\gamma)}.$$

Note that for any density level set estimator $\hat{G}$, it holds that $d_H(\hat{G}, \Gamma(\lambda)) = d_H(\hat{G}, \overline{\Gamma}(\lambda))$. In other words, the choice of definition of the density level set will not affect the performance of an estimator when measured by its excess-mass deficit. However, the distance $d_\Delta$ is very sensitive to this choice, as illustrated in Section 2.2, and one must resort to offsets to control the first term on the right-hand side of the result in Proposition 2.1.

## 2.2. Plug-in density level set estimators with offset

For a fixed $\lambda > 0$, the plug-in estimator of $\Gamma(\lambda)$ is defined by

$$\hat{\Gamma}(\lambda) = \{x \in \mathcal{X} : \hat{p}_n(x) > \lambda\},$$

where $\hat{p}_n$ is a nonparametric estimator of $p$. For example, $\hat{p}_n$ can be a kernel density estimator of $p$,

$$\hat{p}_n(x) = \hat{p}_{n,h}(x) = \frac{1}{nh^d} \sum_{i=1}^{n} K\left(\frac{X_i - x}{h}\right), \qquad x \in \mathcal{X},$$

where $K : \mathbb{R}^d \to \mathbb{R}$ is a suitably chosen kernel and $h > 0$ is the bandwidth parameter. For reasons that will be made clear later, we consider the family of *plug-in estimators with offset* $\ell_n$, denoted by $\tilde{\Gamma}_{\ell_n}$ and defined as

$$\tilde{\Gamma}_{\ell_n} = \tilde{\Gamma}_{\ell_n}(\lambda) = \hat{\Gamma}(\lambda + \ell_n) = \{x \in \mathcal{X} : \hat{p}_n(x) \geq \lambda + \ell_n\},$$

where $\ell_n$ is a quantity that typically tends to 0 as $n$ tends to infinity.

As mentioned in Remark 1.1, when the goal is to estimate the set $\Gamma(\lambda)$, the offset $\ell_n$ is chosen to be positive, whereas for $\overline{\Gamma}$, it must be chosen negative. The effect of such choices is to ensure that the set $\{p = \lambda\}$ is, respectively, removed or added to the standard plug-in estimator with high probability. This phenomenon emerges only when the performance is measured using pseudo-distance $d_\Delta$, but not $d_H$. This translates into different optimal choices for the offset, as well as different optimal rates of convergence depending on the chosen measure of performance.

The following counterexample suggested by an anonymous referee demonstrates that standard plug-in estimators can fail to consistently estimate the set $\{p = \lambda\}$. Assume that $\mathcal{X} \subset \mathbb{R}$, that the density $p$ is such that $p(x) = 1/2$ for all $x \in [0,1]$ and that $p(x) < 1/2$ elsewhere. In this case, it is clear that $\Gamma(1/2) = \varnothing$ and $\overline{\Gamma}(1/2) = [0,1]$. Assume, now, that we have an estimator $\hat{p}$ such that $|\hat{p}(x) - p(x)| \leq \varepsilon$, where $\varepsilon > 0$ is arbitrary small. If $\hat{p}(x) = 1/2 + \varepsilon$ for any $x \in [0,1]$, then $\hat{\Gamma}(1/2) \supset [0,1]$ and it fails to consistently estimate $\Gamma(1/2)$ as $\varepsilon$ tends to 0. However, $\tilde{\Gamma}_{\ell_n}$ with a positive offset $\ell_n > \varepsilon$ can become consistent, as shown in Section 3. Conversely, if $\hat{p}(x) = 1/2 - \varepsilon$ for any $x \in [0,1]$, then $\hat{\Gamma}(1/2)$ is not a consistent estimator of $\overline{\Gamma}(1/2)$, but $\tilde{\Gamma}_{\ell_n}$ with a negative offset $\ell_n < -\varepsilon$ can be one.

As a consequence, plug-in density level set estimators can match both definitions of density level sets (1.1) or (1.2) by simply changing the sign of the offset.

## 3. Fast rates for plug-in density level set estimators with offset

The first theorem states that rates of convergence for plug-in estimators with offset can be obtained using exponential inequalities for the corresponding nonparametric density estimator $\hat{p}_n$. In what follows, smoothness in the neighborhood of the level under consideration is particularly important and we define this neighborhood as follows:

$$\mathcal{D}(\eta) = \{p \in (\lambda - \eta, \lambda + \eta)\}, \qquad \eta > 0.$$

In the sequel, we write, for simplicity, $\tilde{\Gamma}_{\ell_n} = \tilde{\Gamma}$ when the value of the offset is clear from the context.

The following definition will be at the center of our main theorem. It provides a compact way to describe the pointwise convergence of an estimator to the true density $p$.

**Definition 3.1.** *Let $\mathcal{P}$ be a given class of probability densities of $\mathcal{X}$ and fix $\Delta > 0$. Let $\varphi = (\varphi_n)$ and $\psi = (\psi_n)$ be two positive, monotonically non-increasing sequences.*

*We say that an estimator $\hat{p}_n$ is pointwise convergent at a rate $(\psi_n)$ uniformly over $\mathcal{P}$ if there exist positive constants $c_1, c_2, c_\psi$ such that for $Q$-almost all $x \in \mathcal{X}$, we have*

$$\sup_{p \in \mathcal{P}} \mathbb{P}(|\hat{p}_n(x) - p(x)| \geq \delta) \leq c_1 \mathrm{e}^{-c_2 (\delta/\psi_n)^2}, \qquad c_\psi \psi_n < \delta < \Delta. \tag{3.1}$$

*Moreover, we say that an estimator $\hat{p}_n$ is $(\varphi, \psi)$-locally pointwise convergent in $\mathcal{D}(\eta)$, uniformly over $\mathcal{P}$, if it is pointwise convergent at a rate $(\psi_n)$ uniformly over $\mathcal{P}$ and there exists positive constants $c_3, c_4, c_\varphi$ such that for $Q$-almost all $x \in \mathcal{D}(\eta)$, we have*

$$\sup_{p \in \mathcal{P}} \mathbb{P}(|\hat{p}_n(x) - p(x)| \geq \delta) \leq c_3 \mathrm{e}^{-c_4 (\delta/\varphi_n)^2}, \qquad c_\varphi \varphi_n < \delta < \Delta. \tag{3.2}$$

The following theorem states that it is possible to derive fast rates of convergence under the $\gamma$-exponent for plug-in estimators constructed from a locally pointwise convergent estimator of the density.

**Theorem 3.1.** *Fix $\lambda > 0, \Delta > 0$ and let $\mathcal{P}$ be a class of densities on $\mathcal{X}$. Let $\varphi = (\varphi_n)$ and $\psi = (\psi_n)$ be two positive, monotonically non-increasing sequences such that*

$$\varphi_n \vee \psi_n = \mathrm{o}\left(\frac{1}{\sqrt{\log n}}\right) \quad \text{and} \quad \varphi_n \geq C n^{-\mu} \quad \text{for some } C, \mu > 0.$$

*Let $\hat{p}_n$ be an estimator of the density $p$ constructed from data arising from $p \in \mathcal{P}$ and such that $Q(\hat{p}_n \geq \lambda) \leq M$, almost surely for some positive constant $M$. Assume that $\hat{p}_n$ is $(\varphi, \psi)$-locally pointwise convergent in $\mathcal{D}(\eta)$, uniformly over $\mathcal{P}$. Then, if $p$ has $\gamma$-exponent at level $\lambda$ for any $p \in \mathcal{P}$, the plug-in estimator $\tilde{\Gamma}$ based on $\hat{p}_n$, with offset $\ell_n$, satisfies*

$$\sup_{p \in \mathcal{P}} \mathbb{E}[d_H(\Gamma_p(\lambda), \tilde{\Gamma})] \leq C \varphi_n^{(1+\gamma)} \qquad \text{for } \ell_n \leq C \varphi_n, \tag{3.3}$$

$$\sup_{p \in \mathcal{P}} \mathbb{E}[d_\Delta(\Gamma_p(\lambda), \tilde{\Gamma})] \leq C (\varphi_n \sqrt{\log n})^\gamma \qquad \text{for } \ell_n = c_\ell \varphi_n \sqrt{\log n} \tag{3.4}$$

*for $n \geq n_0 = n_0(\lambda, \eta, \varphi, \psi, \varepsilon_0, c_\varphi, c_\psi)$ and where, in (3.4), the constant $c_\ell$ must be chosen large enough so that $c_\ell^2 \geq \mu \gamma / c_4$.*

Before giving the proof of the theorem, we comment on its meaning. First, note that the main consequence of (3.2) is that $|\hat{p}_n(x) - p(x)|$ is of order $\varphi_n$ for any $x$ in the neighborhood $\mathcal{D}(\eta)$. That is, $\hat{p}_n$ is a good pointwise estimator of $p$ in this neighborhood. Equation (3.1) is of the same flavor as (3.2) but in a weaker form. It entails that for $x$ outside $\mathcal{D}(\eta)$, $\hat{p}_n(x)$ is a consistent estimator of $p(x)$ with rate of order $\psi_n$, which can be as slow as $\mathrm{o}((\log n)^{-1/2})$ since it does not appear in the rates (3.3) or (3.4). This confirms the intuition that the density needs to be accurately estimated only in a neighborhood of $\lambda$, whereas, outside this neighborhood, it is sufficient to know whether the density is greater or less that $\lambda$. Finally, note that the constant $c_\ell$ can be constructed from an upper bound on the parameter $\gamma$ – and not $\gamma$ itself – which will be available in the particular context of Section 4. Therefore, $\tilde{\Gamma}$ remains adaptive to the parameter $\gamma$.

As already mentioned in Section 2.2, we can see that the offset employed when the performance is measured using pseudo-distance $d_H$ is negligible with respect the one employed with $d_\Delta$. Actually, a plug-in estimator without offset does the job just as well when $d_H$ is employed.

The proof of the theorem relies on the following lemma. Its proof is postponed to Section A.1 in the Appendix.

**Lemma 3.1.** *Under the assumptions of Theorem 3.1, for any offset $\ell > 0$, the plug-in estimator with offset $\ell > 0$ satisfies*

$$\sup_{p \in \mathcal{P}} \mathbb{E}[d_H(\Gamma_p(\lambda), \tilde{\Gamma}_\ell)] \leq C(\varphi_n \vee \ell)^{(1+\gamma)}, \tag{3.5}$$

*for some positive constant $C$.*

We now turn to the proof of Theorem 3.1.

**Proof of Theorem 3.1.** Note first that (3.3) is a direct consequence of Lemma 3.1 applied with $\ell = \ell_n \leq C\varphi_n$.

To prove (3.4), we apply Proposition 2.1, whose conditions are satisfied. Indeed, using the Markov inequality, we get

$$Q(\{\hat{p}_n \geq \lambda + \ell_n\}\Delta\{p \geq \lambda\}) \leq Q(\hat{p}_n \geq \lambda) + Q(p \geq \lambda) \leq M + \lambda^{-1}$$

and we choose $L_Q = M + \lambda^{-1}$. Proposition 2.1 yields

$$\mathbb{E}[d_\Delta(\Gamma_p(\lambda), \tilde{\Gamma})] \leq \mathbb{E}Q(\Gamma_p(\lambda)\Delta\tilde{\Gamma} \cap \{p = \lambda\}) + C(\varphi_n\sqrt{\log n})^\gamma, \qquad (3.6)$$

where the second term in the right-hand side was controlled by successively applying the Jensen inequality and Lemma 3.1 with $\ell = c_\ell \varphi_n \sqrt{\log n}$, as prescribed in (3.4). To conclude the proof, it is sufficient to observe that by Fubini's Theorem and assumption (3.2), we have

$$\mathbb{E}Q(\Gamma_p(\lambda)\Delta\tilde{\Gamma} \cap \{p = \lambda\}) \leq c_3 Q(p = \lambda)\mathrm{e}^{-c_4(\ell_n/\varphi_n)^2} \leq Cn^{-\mu\gamma} \leq C\varphi_n^\gamma.$$

Together with (3.6), this inequality yields (3.4), which concludes the proof. □

In the next section, we verify that kernel density estimators are $(\varphi, \psi)$-locally pointwise convergent in a neighborhood of $\lambda$, uniformly over a class of locally Hölder smooth probability densities.

## 4. Optimal rates for plug-in estimators with offset based on kernel density estimators

In the remainder of this paper, we fix the measure $Q$ to be the Lebesgue measure on $\mathbb{R}^d$, denoted by $\mathrm{Leb}_d$.

In this section, we derive exponential inequalities of type (3.1) when the estimator $\hat{p}_n$ is a kernel density estimator and the density $p$ belongs to some Hölder class of densities. We begin by giving the definition of the Hölder classes of densities that we consider.

### 4.1. Hölder classes of densities

Fix $\beta > 0$ and $\lambda > 0$. For any $d$-tuples $s = (s_1, \ldots, s_d) \in \mathbb{N}^d$ and $x = (x_1, \ldots, x_d) \in \mathcal{X}$, we define $|s| = s_1 + \cdots + s_d$, $s! = s_1! \cdots s_d!$ and $x^s = x_1^{s_1} \cdots x_d^{s_d}$. Let $D^s$ denote the differential operator

$$D^s = \frac{\partial^{s_1 + \cdots + s_d}}{\partial x_1^{s_1} \cdots \partial x_d^{s_d}}.$$

For any real-valued function $g$ on $\mathcal{X}$ that is $\lfloor \beta \rfloor$-times continuously differentiable at point $x_0 \in \mathcal{X}$, we denote by $g_{x_0}^{(\beta)}$ its Taylor polynomial of degree $\lfloor \beta \rfloor$ at point $x_0$:

$$g_{x_0}^{(\beta)}(x) = \sum_{|s| \leq \lfloor \beta \rfloor} \frac{(x-x_0)^s}{s!} D^s g(x_0).$$

Fix $L > 0, r > 0$ and denote by $\Sigma(\beta, L, r, x_0)$ the set of functions $g: \mathcal{X} \to \mathbb{R}$ that are $\lfloor \beta \rfloor$-times continuously differentiable at point $x_0$ and satisfy

$$|g(x) - g_{x_0}^{(\beta)}(x)| \leq L\|x - x_0\|^\beta \qquad \forall x \in \mathcal{B}(x_0, r).$$

The set $\Sigma(\beta, L, r, x_0)$ is called the $(\beta, L, r, x_0)$-*locally Hölder class* of functions. We now define the class of densities that are considered in this paper.

**Definition 4.1.** *Fix $\beta > 0, L > 0, r > 0, \lambda > 0$ and $\gamma > 0$. Recall that $\mathcal{D}(\eta)$ is the neighborhood defined by*

$$\mathcal{D}(\eta) = \{p \in (\lambda - \eta, \lambda + \eta)\}, \qquad \eta > 0.$$

*Let $\mathcal{P}_\Sigma(\beta, L, r, \lambda, \gamma, \beta', L^*)$ denote the class of all probability densities $p$ on $\mathcal{X}$ for which there exists $\eta > 0$ such that:*

(i) *$p \in \Sigma(\beta, L, r, x_0)$ for all $x_0 \in \mathcal{D}(\eta)$, apart from a set of null Lebesgue measure;*
(ii) *$\exists \beta' > 0$ such that $p \in \Sigma(\beta', L, r, x_0)$ for all $x_0 \notin \mathcal{D}(\eta)$, apart from a set of null Lebesgue measure;*
(iii) *$p$ has $\gamma$-exponent at level $\lambda$ with respect to the Lebesgue measure;*
(iv) *$p$ is uniformly bounded by a constant $L^*$.*

*We will often prefer the compact notation $\mathcal{P}_\Sigma(\beta, \lambda, \gamma)$, or simply $\mathcal{P}_\Sigma$, when either the parameters are clear from the context or their value does not affect the results.*

The class $\mathcal{P}_\Sigma(\beta, \lambda, \gamma)$ is the class of uniformly bounded (iv) densities that have $\gamma$-exponent at level $\lambda$ with respect to $\text{Leb}_d$ (iii) and that are smooth in the neighborhood of the level under consideration (i). As usual in nonparametric estimation, the parameters $L, L^*$ and $r$ will affect only the constants in the rates of convergence presented below. However, the smoothness parameter $\beta'$ in condition (ii) which is expected to control the rate of convergence will also affect only the constants and therefore does not appear in the compact notation $\mathcal{P}_\Sigma(\beta, \lambda, \gamma)$. Indeed, $\beta' > 0$ can be arbitrarily close to 0 and this will not affect the rates of convergence. Actually, the role of condition (ii) is to ensure that any density from the class can be consistently estimated at any point with an arbitrarily slow polynomial rate.

The class of densities $\mathcal{P}_\Sigma$ is similar to the class of regression functions considered in Audibert and Tsybakov (2007). However, besides the additional assumption that functions in $\mathcal{P}_\Sigma$ are probability densities, the main improvement here is that the regularity of a density in $\mathcal{P}_\Sigma$ can be arbitrarily low outside a neighborhood of the level under consideration, yielding slower rates of pointwise estimation. We prove below (cf. Corollary 4.1) that optimal rates of convergence for DLSE are possible for this larger class of densities, which corroborates the idea that the density need not be precisely estimated far from the level $\lambda$.

The next proposition can be derived by following the lines of the proof of Proposition 3.4 (fourth item) of Audibert and Tsybakov (2005).

**Proposition 4.1.** *If $\gamma(\beta \wedge 1) > 1$, either $\Gamma$ has empty interior or its complement $\Gamma^c$ does. Conversely, if $\gamma(\beta \wedge 1) \leq 1$, then there exist densities such that both $\Gamma$ and $\Gamma^c$ have non-empty interior.*

### 4.2. Exponential inequalities for kernel density estimators

To estimate a density $p$ from the class $\mathcal{P}_\Sigma(\beta, \lambda, \gamma)$, we can use a *kernel density estimator* defined by

$$\hat{p}_n(x) = \hat{p}_{n,h}(x) = \frac{1}{nh^d} \sum_{i=1}^n K\left(\frac{X_i - x}{h}\right), \qquad (4.1)$$

where $h > 0$ is the bandwidth parameter and $K : \mathcal{X} \to \mathbb{R}$ is a kernel. This choice is not the only possible one and all we need is an estimator that satisfies exponential inequalities as in (3.1) and (3.2). The following lemma states that it is possible to derive such exponential inequalities for a kernel density estimator with a $\beta^\star$-valid kernel where $\beta^\star \geq \beta$. The definition of $\beta$-valid kernel is recalled in the Appendix, Definition A.1 (see also Tsybakov (2004) for example).

**Lemma 4.1.** *Let $P$ be a distribution on $\mathbb{R}^d$ having a density $p$ with respect to the Lebesgue measure and such that $\|p\|_\infty \leq L^*$ for some constant $L^* > 0$. Fix $\beta > 0$, $\beta^\star \geq \beta$, $L > 0$, $r > 0$ and assume that $p \in \Sigma(\beta, L, r, x_0)$. Let $\hat{p}_n$ be a kernel density estimator with bandwidth $h > 0$ and $\beta^\star$-valid kernel $K$, given an i.i.d. sample $X_1, \ldots, X_n$ from $P$. Set*

$$\Delta = \frac{6L^* \|K\|^2}{\|K\|_\infty + L^* + L \int \|t\|^\beta K(t) \, \mathrm{d}t}.$$

*Then, for all $\delta, h \leq r$ such that $\Delta > \delta > 2Lc_5 h^\beta > 0$, we have*

$$\mathbb{P}\{|\hat{p}_n(x_0) - p(x_0)| \geq \delta\} \leq 2 \exp(-c_6 n h^d \delta^2),$$

*where $c_5 = \int \|t\|^\beta K(t) \, \mathrm{d}t$ and $c_6 = 1/(16L^*\|K\|^2)$.*

The proof is given in Section A.2 of the Appendix.

We can therefore apply Theorem 3.1. For the appropriate choice of $h$, it yields the following corollary.

**Corollary 4.1.** *Let $Q$ be the Lebesgue measure on $\mathbb{R}^d$. Fix positive constant $\beta, L, r, \lambda, \gamma, \beta', L^*$ and assume that $\gamma(\beta \wedge 1) \leq 1$. Consider the class of densities $\mathcal{P}_\Sigma = \mathcal{P}_\Sigma(\beta, L, r, \lambda, \gamma, \beta', L^*)$ and define $\beta^\star = \beta \vee \beta'$.*

Let $\tilde{\Gamma}$ be the plug-in estimator with offset $\ell_n$ based on the estimator $\hat{p}_n$ defined in (4.1) with $\beta^\star$-valid kernel $K$ and bandwidth parameter $h_n > 0$. Then, for any $c_h > 0$ and for $c_\ell \geq \max((c_6 c_h^d)^{-1}, 1)$, we have

$$\sup_{p \in \mathcal{P}_\Sigma} \mathbb{E}[d_H(\Gamma_p(\lambda), \tilde{\Gamma})] \leq C n^{-(1+\gamma)\beta/(2\beta+d)} \quad \text{for } \begin{cases} \ell_n \leq C n^{-\beta/(2\beta+d)}, \\ h_n = c_h n^{-1/(2\beta+d)}; \end{cases}$$

$$\sup_{p \in \mathcal{P}_\Sigma} \mathbb{E}[d_\Delta(\Gamma_p(\lambda), \tilde{\Gamma})] \leq C \left(\frac{n}{\log n}\right)^{-\gamma\beta/(2\beta+d)} \quad \text{for } \begin{cases} \ell_n = c_\ell n^{-\beta/(2\beta+d)} \sqrt{\log n}, \\ h_n = c_h (n/\log n)^{-1/(2\beta+d)}. \end{cases}$$

**Proof.** Define $\varphi_n = (n h_n^d)^{-1/2}$ and $\psi_n = h_n^{\beta'} \geq C h_n^\beta$ and consider separately the cases $h_n = c_h n^{-1/(2\beta+d)}$ and $h_n = c_h (n/\log n)^{-1/(2\beta+d)}$. We have, respectively,

$$\varphi_n = c_\varphi n^{-\beta/(2\beta+d)} \leq C h_n^\beta \quad \text{and} \quad \varphi_n = C n^{-\beta/(2\beta+d)} (\log n)^{-(d/2)/(2\beta+d)} < C h_n^\beta.$$

Therefore, in both cases, $\psi_n \geq C h_n^\beta \geq C(n h_n^d)^{-1/2}$ and a direct consequence of Lemma 4.1 is that the kernel density estimator with $\beta^\star$-valid kernel $K$ and bandwidth parameter $h_n > 0$ is $(\varphi, \psi)$-locally pointwise convergent in $\mathcal{D}(\eta)$, uniformly over $\mathcal{P}_\Sigma$ in both cases.

We also need to check that for such an estimator we have $\text{Leb}_d(\hat{p}_n \geq \lambda) \leq M$, almost surely for some $M > 0$. Note that since $K \in L_1(\mathbb{R}^d)$, we have

$$\infty > \int_{\mathbb{R}^d} |K(x)| \, d\text{Leb}_d(x) \geq \int_{\{\hat{p}_n \geq \lambda\}} |\hat{p}_n(x)| \, d\text{Leb}_d(x) \geq \lambda \, \text{Leb}_d\{\hat{p}_n \geq \lambda\}.$$

Hence, the condition is satisfied with $M = \lambda^{-1} \int |K|$.

Let us finally check that the aforementioned choice of $c_\ell$ is compatible with the assumptions of Theorem 3.1. Since $\varphi_n = c_\varphi n^{-\beta/(2\beta+d)}$, we can take $\mu = \beta/(2\beta+d)$. We need to check that $c_\ell^2 \geq \mu \gamma/(c_6 c_h^d)$. Since $c_\ell \geq 1$, we have $c_\ell^2 \geq c_\ell$ and

$$\frac{\mu \gamma}{c_6 c_h^d} \leq \frac{\beta}{(2\beta+d)(\beta \wedge 1) c_6 c_h^d} = \frac{\beta \vee 1}{(2\beta+d)(c_6 c_h^d)} \leq \frac{1}{c_6 c_h^d} \leq c_\ell.$$

Since all the conditions of Theorem 3.1 are satisfied, this concludes the proof. □

## 5. Minimax lower bounds

The following theorem shows that the rates obtained in Corollary 4.1 are optimal in a minimax sense.

**Theorem 5.1.** *Let the underlying measure $Q$ be the Lebesgue measure on $\mathbb{R}^d$. Fix $\lambda > 0$, let $\beta, \gamma$ be positive constants such that $\gamma \beta \leq 1$ and consider the class of densities $\mathcal{P}_\Sigma = \mathcal{P}_\Sigma(\beta, \lambda, \gamma)$. Then, for any $n \geq 1$ and any estimator $\hat{G}_n$ of $\Gamma_p(\lambda)$ constructed from the* sample $X_1, \ldots, X_n$, we have

$$\sup_{p \in \mathcal{P}_\Sigma} \mathbb{E}[d_H(\Gamma_p(\lambda), \hat{G}_n)] \geq C n^{-(1+\gamma)\beta/(2\beta+d)}, \tag{5.1}$$

$$\sup_{p \in \mathcal{P}_\Sigma} \mathbb{E}[d_\Delta(\Gamma_p(\lambda), \hat{G}_n)] \geq C \left(\frac{n}{\log n}\right)^{-\gamma\beta/(2\beta+d)}. \tag{5.2}$$

**Proof.** Fix $n \geq 2$ and consider the quantities

$$\varrho_n = C n^{-\gamma\beta/(2\beta+d)} \quad \text{and} \quad \varkappa_n = C \left(\frac{n}{\log n}\right)^{-\gamma\beta/(2\beta+d)}.$$

Our goal is to find two families of densities $\mathcal{N}_\varrho$ and $\mathcal{N}_\varkappa$ that are in $\mathcal{P}$ and which satisfy the conditions of Lemma A.2. Note first that although it does not appear in its notation, the pseudo-distance $d_H$ depends on the underlying density $p$ and this will be inconvenient in proving minimax lower bounds. As a result, we are going to prove (5.1) for a pseudo-distance $d_\varrho$ which does not depend on $p$ and, for any measurable $G \subset [0,1]^d$, satisfies

$$d_H(\Gamma_p(\lambda), G) \geq C(d_\varrho(\Gamma_p(\lambda), G))^{(1+\gamma)/\gamma} \qquad \text{for any } p \in \mathcal{N}_\varrho. \tag{5.3}$$

The pseudo-distance $d_\varrho$ will be defined in (5.6). We will use Lemma A.2 for $\mathcal{P} = \mathcal{P}_\Sigma$, with $\varepsilon = \varrho_n$, $\mathcal{N} = \mathcal{N}_\varrho$, $d = d_\varrho$ to prove (5.1), and with $\varepsilon = \varkappa_n$, $\mathcal{N} = \mathcal{N}_\varkappa$, $d = d_\Delta$ to prove (5.2).

For both families $\mathcal{N}_\varrho$ and $\mathcal{N}_\varkappa$, the construction begins as follows. Assume, without loss of generality, that $\lambda = 1$. Let $q \geq 4$ be an integer, to be specified later, and consider the regular grid $\mathcal{G}$ on $[0,1]^d$ defined as

$$\mathcal{G} = \left\{\left(\frac{2k_1+1}{q}, \ldots, \frac{2k_d+1}{q}\right), k_i \in \{0, \ldots, q-1\}, i = 1, \ldots, d\right\}.$$

Let $N$ denote the unique integer in the pair $\{q^d/2, (q^d-1)/2\}$ and denote by $\{g_j\}_{1 \leq j \leq 2N}$ a collection of $2N$ distinct elements of the grid, the choice of indexing being of no importance for what follows. For any $j = 1, \ldots, 2N$, define the Euclidean balls $B_j = \mathcal{B}(g_j, \kappa)$, where $\kappa = 1/q$.

Let $\phi_\beta : \mathbb{R}^d \to \mathbb{R}_+$ be a smooth function defined as follows. If $\beta < 1$, the function $\phi_\beta$ is defined as

$$\phi_\beta(x) = \begin{cases} C_\beta(1 - \|x\|)^\beta, & \text{if } 0 \leq \|x\| \leq 1, \\ 0, & \text{if } \|x\| > 1. \end{cases}$$

If $\beta \geq 1$, the function $\phi_\beta$ is defined as

$$\phi_\beta(x) = \begin{cases} C_\beta(2^{1-\beta} - \|x\|^\beta), & \text{if } 0 \leq \|x\| \leq 1/2, \\ C_\beta(1 - x)^\beta, & \text{if } 1/2 \leq \|x\| \leq 1, \\ 0, & \text{if } \|x\| > 1, \end{cases}$$

where, in both cases, $0 < C_\beta < 1/2$ is chosen small enough to ensure that $|\phi_\beta(x) - \phi_\beta(x')| \leq L\|x - x'\|^\beta$ for any $x, x' \in \mathbb{R}^d$.

Then, for any $\omega = (\omega_1, \ldots, \omega_N) \in \{-1, 0, 1\}^N$, define on $[0,1]^d$ the function

$$p_\omega(x) = 1 + \sum_{j=1}^{N} \omega_j [\varphi_j(x) - \varphi_{N+j}(x)],$$

where $\varphi_j(x) = \kappa^\beta \phi([x - g_j]/\kappa) \mathbb{1}_{\{x \in B_j\}}$.

Define the integer $m = \lfloor q^{d-\gamma\beta}/2 \rfloor + 6$ so that $m$ satisfies $6 \leq m \leq 6N$.

Let $\Omega_\varrho$ and $\Omega_\varkappa$ be two subsets of $\{-1, 0, 1\}^N$ such that $\sum_j |\omega_j| = m$ for any $\omega = (\omega_1, \ldots, \omega_N) \in \Omega_\varrho \cup \Omega_\varkappa$. Now define the families $\mathcal{N}_\varrho$ and $\mathcal{N}_\varkappa$ as

$$\mathcal{N}_\varrho = \{p_\omega, \omega \in \Omega_\varrho\}, \qquad \mathcal{N}_\varkappa = \{p_\omega, \omega \in \Omega\varkappa\}.$$

The sets $\Omega_\varrho$ and $\Omega_\varkappa$, will be chosen in order to fulfill the conditions of Lemma A.2.

***First condition.*** $\mathcal{N} \subset \mathcal{P}_\Sigma(\beta, 1, \gamma)$.

First, note that for any $\omega \in \{-1, 0, 1\}^N$, $\|p_\omega\|_\infty \leq 2$ and $p_\omega \in \Sigma(\beta, 2, 1, x)$ for any $x \in [0,1]^d$. Therefore, it remains to check that $p_\omega$ has $\gamma$-exponent at level 1 with respect to the Lebesgue measure. We now show that it is sufficient to have $\sum_j |\omega_j| \leq 2m$ for this condition to hold, which is satisfied for $p_\omega$ either in $\mathcal{N}_\varrho$ or in $\mathcal{N}_\varkappa$. We have

$$\mathrm{Leb}_d(x: 0 < |p_\omega(x) - 1| \leq \varepsilon) = 2 \sum_{j=1}^{N} \mathbb{1}_{\{|\omega_j|=1\}} \mathrm{Leb}_d(x: 0 < |p_\omega(x) - 1| \leq \varepsilon, x \in B_j)$$

$$\leq 4m \, \mathrm{Leb}_d(x: 0 < \phi([x - g_1]/\kappa) \leq \varepsilon \kappa^{-\beta})$$

$$= 4m \int_{\mathcal{B}(0,1)} \mathbb{1}_{\{\phi(x/\kappa) \leq \varepsilon \kappa^{-\beta}\}} \, dx.$$

The last term in the previous system of equations is treated differently depending on whether $\beta < 1$ or $\beta \geq 1$. If $\beta < 1$, we have

$$\begin{aligned}
4m \int_{\mathcal{B}(0,1)} \mathbb{1}_{\{\phi(x/\kappa) \leq \varepsilon \kappa^{-\beta}\}} \, dx &\leq C(m\kappa^d \mathbb{1}_{\{\varepsilon > \kappa^\beta\}} + m[\kappa^d - (\kappa - \varepsilon^{1/\beta})^d] \mathbb{1}_{\{\varepsilon \leq \kappa^\beta\}}) \\
&\leq C(\kappa^{\gamma\beta} \mathbb{1}_{\{\varepsilon^\gamma > \kappa^{\gamma\beta}\}} + m\kappa^{d-1} \varepsilon^{1/\beta} \mathbb{1}_{\{\kappa > \varepsilon^{1/\beta}\}}) \\
&\leq C(\varepsilon^\gamma + \kappa^{\gamma\beta - 1} \varepsilon^{1/\beta} \mathbb{1}_{\{\kappa > \varepsilon^{1/\beta}\}}) \\
&\leq C\varepsilon^\gamma,
\end{aligned} \qquad (5.4)$$

where we have used the fact that $\gamma\beta - 1 \leq 0$ to bound the second term in the penultimate inequality.

We now treat the case $\beta > 1$. Note that integration over $x$ such that $\|x\| \geq \kappa/2$ can be treated in the same manner as in the case $\beta < 1$. The integral over $x$ such that $\|x\| < \kappa/2$ is trivially upper bounded by a term proportional to the volume of the ball $\mathcal{B}(0, \kappa/2)$. It yields

$$4m \int_{\mathcal{B}(0,1)} \mathbb{1}_{\{\phi(x/\kappa) \leq \varepsilon \kappa^{-\beta}\}} \, \mathrm{d}x \leq C(m\kappa^d \mathbb{1}_{\{\varepsilon > C_\beta(\kappa/2)^\beta\}} + m\kappa^{d-1}\varepsilon^{1/\beta} \mathbb{1}_{\{\varepsilon \leq C_\beta(\kappa/2)^\beta\}})$$
$$\leq C(\kappa^{\gamma\beta} \mathbb{1}_{\{\kappa \leq 2(\varepsilon/C_\beta)^{1/\beta}\}} + m\kappa^{d-1}\varepsilon^{1/\beta} \mathbb{1}_{\{\kappa > 2(\varepsilon/C_\beta)^{1/\beta}\}}) \quad (5.5)$$
$$\leq C\varepsilon^\gamma.$$

As a result, both $\mathcal{N}_\varrho$ and $\mathcal{N}_\varkappa$ are subsets of $\mathcal{P}_\Sigma(\beta, 1, \gamma)$ and the first condition is satisfied.

**Second condition.** $d(\Gamma_p, \Gamma_q) \geq \varepsilon \, \forall \, p, q \in \mathcal{N}, p \neq q$.

It will become clearer that we need to bound from below the Hamming distance between $\omega$ and $\omega'$, defined for any $\omega, \omega' \in \{-1, 0, 1\}^N$ by

$$\rho(\omega, \omega') = \sum_{j=1}^{N} \mathbb{1}_{\{\omega_j \neq \omega'_j\}}.$$

Let us now treat separately the class $\mathcal{N}_\varrho$ and the class $\mathcal{N}_\varkappa$.

We begin with $\mathcal{N}_\varrho$ and define $\tilde{l} = \lfloor m/6 \rfloor$ so that $m \geq 6\tilde{l} \geq 6$. The subset $\Omega_\varrho$ is chosen to be of the form

$$\Omega_\varrho = \{\omega = (\omega^{[m]}, 0, \ldots, 0), \omega^{[m]} \in \Omega_\varrho^{[m]}\},$$

where $\Omega_\varrho^{[m]}$ is a subset of $\{-1, 1\}^m$. For any $\omega \in \Omega_\varrho$, we clearly have

$$\sum_{j=1}^{N} |\omega_j| = m \leq 2m.$$

We are now in a position to define $d_\varrho$ by

$$d_\varrho(G_1, G_2) = \sum_{j=1}^{m} \mathrm{Leb}_d(G_1 \Delta G_2 \cap (B_j \cup B_{N+j})). \quad (5.6)$$

It is easy to check that (5.3) is satisfied since for any $p \in \mathcal{N}_\varrho$ and any measurable $G \subset [0,1]^d$, we have

$$d_\varrho(\Gamma_p, G) = \mathrm{Leb}_d(\Gamma_p \Delta G \cap \{p \neq \lambda\}) \leq C(d_H(\Gamma_p(\lambda), G))^{\gamma/(1+\gamma)},$$

where in the last inequality, we have used Proposition 2.1.

The set $\Omega_\varrho^{[m]}$ is extracted from $\{-1,1\}^m$ using Lemma A.1, which guarantees that there exists such a set with cardinality $s_\varrho \geq 2$ satisfying

$$\log(s_\varrho) \geq C\tilde{l}\log(m/\tilde{l}) \geq Cm$$

and

$$\sum_{j=1}^m \mathbb{1}_{\{\omega_j \neq \omega_j'\}} \geq \tilde{l}+1 \geq m/6.$$

It yields

$$d_\varrho(\Gamma_{p_\omega}, \Gamma_{p_{\omega'}}) = 2\operatorname{Leb}_d(B_1)\rho(\omega,\omega') \geq Cm\kappa^d \geq 2\varrho_n$$

when $q = 4\lfloor n^{1/(2\beta+d)}\rfloor$.

We now define the set $\mathcal{N}_\varkappa$ as follows. Let $\Omega_\varkappa$ be a subset of $\{0,1\}^N$ with cardinality $s_\varkappa \geq 2$, extracted using Lemma A.1 such that

$$\sum_{j=1}^N \omega_j = 2m$$

for any $\omega$ in this extracted subset and to which $0 \in \{0,1\}^N$ has been added. It satisfies

$$\log(s_\varkappa) \geq Cm\log(N/m) \geq Cm\log(q)$$

and

$$\sum_{j=1}^N \mathbb{1}_{\{\omega_j \neq \omega_j'\}} \geq m+1.$$

It yields

$$d_\Delta(\Gamma_{p_\omega}, \Gamma_{p_{\omega'}}) = 2\operatorname{Leb}_d(B_1)\sum_{j=1}^N \mathbb{1}_{\{\omega_j \neq \omega_j'\}} \geq C\kappa^d m \geq 2\varkappa_n$$

when $q = 4\lfloor (n/\log n)^{1/(2\beta+d)}\rfloor$.

**Third condition.** $\max_{\omega \in \Omega} K(p_\omega, p_{\omega_0}) \leq C\log(\operatorname{card}(\mathcal{N}))$.

Note that for the above choice of $\Omega_\varrho$, we have $\log(\operatorname{card}(\mathcal{N}_\varrho)) \geq Cm$ and need only prove that

$$\max_{\omega,\omega' \in \Omega_\varrho} K(p_\omega, p_{\omega'}) \leq Cm.$$

Define $\xi_j(x) = \varphi_j(x) - \varphi_{j+N}(x)$. For any $p_\omega, p_{\omega'} \in \mathcal{N}$, using the inequality

$$\log\left(\frac{1+a}{1+b}\right)(1+a) \leq (a-b) + 2(a-b)^2$$

for any $a, b$ such that $|a| < 1, |b| < 1/2$, we have, for $q = 4\lfloor n^{1/(2\beta+d)} \rfloor$,

$$K(p_\omega, p_{\omega'}) = n \sum_{j=1}^{m} \int_{B_j \cup B_{N+j}} \log\left(\frac{1+\omega_j \xi_j(x)}{1+\omega'_j \xi_j(x)}\right)(1+\omega_j \xi_j(x))\,dx$$

$$\leq n \sum_{j=1}^{m} \int_{B_j \cup B_{N+j}} 2[(\omega_j - \omega'_j)\xi_j(x)]^2\,dx$$

$$\leq 2nm \int_{B_1} \varphi_1^2(x)\,dx$$

$$\leq 2nm\kappa^{(2\beta+d)} \int_{\mathcal{B}(0,1)} \phi^2(x)\,dx$$

$$\leq Cm,$$

where in the third line, we used the invariance by translation of the family $\{\varphi_j\}_j$.

We can therefore apply Lemma A.2 to prove that

$$\sup_{p \in \mathcal{N}_\varrho} \mathbb{E}[d_\varrho(\Gamma_p(\lambda), \hat{G}_n)] \geq C\varrho_n.$$

This inequality combined with (5.3) and the Jensen inequality yields (5.1).

In the same manner, to prove (5.2), it is sufficient to prove that

$$\max_{\omega \in \Omega_\varkappa} K(p_\omega, p_0) \leq Cm \log q.$$

This follows from the following sequence of inequalities:

$$K(p_\omega, p_0) = n \sum_{j=1}^{N} \int_{B_j \cup B_{N+j}} \log(1+\omega_j \xi_j(x))(1+\omega_j \xi_j(x))\,dx$$

$$\leq 2n \sum_{j=1}^{N} \int_{B_1} \omega_j^2 \varphi_1^2(x)\,dx$$

$$\leq 2nm\kappa^{(2\beta+d)} \int_{\mathcal{B}(0,1)} \phi^2(x)\,dx$$

$$\leq Cm \log n$$

$$\leq Cm \log q$$

for $q = 4\lfloor n^{1/(2\beta+d)} \rfloor$, $m = C(n/\log n)^{(d-\gamma\beta)/(2\beta+d)}$ and where, in the first inequality, we have used the convexity inequality $\log(1+x) \leq x$, together with the invariance by translation of the family $\{\varphi_j\}_j$. □

# Appendix

Several results that can be omitted in a first reading of the paper are collected in this appendix.

## A.1. Proof of Lemma 3.1

To prove (3.5), we use the same scheme as in the proof of Audibert and Tsybakov (2007), Theorem 3.1. Recall that $\tilde{\Gamma}$ denotes the plug-in estimator with offset $0 \leq \ell \leq c_\ell \varphi_n \sqrt{\log n}$ and that $\tilde{\Gamma} \Delta \Gamma = (\tilde{\Gamma} \cap \Gamma^c) \cup (\tilde{\Gamma}^c \cap \Gamma)$. It yields

$$\mathbb{E}[d_H(\Gamma, \tilde{\Gamma})] = \mathbb{E} \int_{\tilde{\Gamma} \cap \Gamma^c} |p(x) - \lambda|\, dQ(x) + \mathbb{E} \int_{\tilde{\Gamma}^c \cap \Gamma} |p(x) - \lambda|\, dQ(x).$$

Define two sequences

$$\alpha_n = c_\alpha(\varphi_n \vee \ell) \quad \text{and} \quad \beta_n = c_\beta(\varphi_n \vee \psi_n)\sqrt{\log n}, \qquad n \geq 2,$$

where $c_\alpha = 2\max(c_\varphi, c_\psi, 1)$ and $c_\beta \geq c_\alpha \max(c_\ell, 2\mu(1+\gamma)/c_2, 1)$. Let $n_0$ be a positive integer such that $\alpha_n < \beta_n < \eta \wedge \varepsilon_0 \wedge \Delta$ for all $n \geq n_0$. In the remainder of the proof, we always assume that $n \geq n_0$. Consider the following disjoint decomposition:

$$\tilde{\Gamma}^c \cap \Gamma = \{\hat{p}_n < \lambda + \ell, p > \lambda\} \subset A_1 \cup A_2 \cup A_3, \tag{A.1}$$

where

$$A_1 = \{\hat{p}_n < \lambda + \ell, \lambda < p \leq \lambda + \alpha_n\},$$
$$A_2 = \{\hat{p}_n < \lambda + \ell, \lambda + \alpha_n < p \leq \lambda + \beta_n\},$$
$$A_3 = \{\hat{p}_n < \lambda + \ell, p > \lambda + \beta_n\}.$$

Observe that $A_1 \subseteq \{0 < |p - \lambda| \leq \alpha_n\}$. This yields

$$\mathbb{E} \int_{A_1} |p(x) - \lambda|\, dQ(x) \leq \alpha_n Q(A_1) \leq c_0(\alpha_n)^{1+\gamma}, \tag{A.2}$$

where, in the last inequality, we used the $\gamma$-exponent of $p$. Define $J_n = \lfloor \log_2(\frac{\beta_n}{\alpha_n}) \rfloor + 2$, where $\lfloor y \rfloor$ denotes the maximal integer that is strictly smaller than $y > 0$. We can then partition $A_2$ into

$$A_2 = \bigcup_{j=1}^{J_n} \mathcal{X}_j \cap A_2,$$

where

$$\mathcal{X}_j = \{\hat{p}_n < \lambda + \ell, \lambda + 2^{j-1}\alpha_n < p \leq \lambda + 2^j \alpha_n\} \cap \mathcal{D}(\eta \wedge \varepsilon_0).$$

Hence,

$$\mathbb{E}\int_{A_2}|p(x)-\lambda|\,\mathrm{d}Q(x) = \sum_{j=1}^{J_n}\mathbb{E}\int_{\mathcal{X}_j\cap A_2}|p(x)-\lambda|\,\mathrm{d}Q(x). \tag{A.3}$$

Now, since $\ell \leq \alpha_n/2$, we have

$$\mathcal{X}_j \subset \{|\hat{p}_n - p| > 2^{j-2}\alpha_n\}\cap\{|p(x)-\lambda| < 2^j\alpha_n\}.$$

Using Fubini's theorem and the previous inclusion, the general term of the sum in the right-hand side of (A.3) can be bounded from above by

$$2^j\alpha_n\int_{\mathcal{D}(\eta\wedge\varepsilon_0)}\mathbb{P}[|\hat{p}_n(x)-p(x)|>2^{j-2}\alpha_n]\mathbb{1}_{\{0<|p(x)-\lambda|<2^j\alpha_n\}}\,\mathrm{d}Q(x).$$

Note that for any $1 \leq j \leq J_n$, we have $c_\varphi\varphi_n \leq 2^{j-2}\alpha_n \leq \beta_n \leq \Delta$. Now using (3.2) and the fact that $p$ has $\gamma$-exponent at level $\lambda$, we get

$$\begin{aligned}\mathbb{E}\int_{A_2}|p(x)-\lambda|\,\mathrm{d}Q(x) &\leq c_0c_3\sum_{j\geq 1}(2^j\alpha_n)^{1+\gamma}\exp(-c_4(2^{j-2}\alpha_n/\varphi_n)^2)\\ &\leq C(\alpha_n)^{1+\gamma},\end{aligned} \tag{A.4}$$

where we have used the fact that $\varphi_n \leq \alpha_n$.

We now treat the integral over $A_3$ by first noting that $Q(A_3) \leq 1/\lambda$, from the Markov inequality. Next, using Fubini's theorem and the fact that $\beta_n/2 \geq \ell$ and $\beta_n/2 \geq c_\psi\psi_n$, we obtain

$$\begin{aligned}\mathbb{E}\int_{A_3}|p(x)-\lambda|\,\mathrm{d}Q(x) &\leq \int_{A_3}|p(x)-\lambda|\mathbb{P}[|\hat{p}_n(x)-p(x)|>\beta_n/2]\,\mathrm{d}Q(x)\\ &\leq 2c_1\exp(-c_2(\beta_n/(2\psi_n))^2) \leq 2c_1 n^{-c_2 c_\beta^2/2}.\end{aligned}$$

Using the fact that $c_\beta^2 > c_\beta > 2\mu(1+\gamma)/c_2$, we get

$$\mathbb{E}\int_{A_3}|p(x)-\lambda|\,\mathrm{d}Q(x) \leq 2c_1 n^{-\mu(1+\gamma)} \leq C\alpha_n^{(1+\gamma)}. \tag{A.5}$$

In view of (A.1), if we combine (A.2), (A.4) and (A.5), we obtain

$$\mathbb{E}\int_{\tilde{\Gamma}^c\cap\Gamma}|p(x)-\lambda|\,\mathrm{d}Q(x) \leq C\alpha_n^{(1+\gamma)}.$$

In the same manner, it can be shown that for $n \geq n_0$,

$$\mathbb{E}\int_{\tilde{\Gamma}\cap\Gamma^c}|p(x)-\lambda|\,\mathrm{d}Q(x) \leq C\alpha_n^{(1+\gamma)}.$$

The only difference with the part of the proof detailed above is that in the step that corresponds to proving the equivalent of (A.5), we use the assumption that $Q(\hat{p}_n \geq \lambda) \leq M$ almost surely, in place of the Markov inequality.

## A.2. Proof of Lemma 4.1

**Proof.** For any $x_0 \in \mathbb{R}^d$,

$$|\hat{p}_n(x_0) - p(x_0)| = \frac{1}{n}\left|\sum_{i=1}^{n} Z_i(x_0)\right|,$$

with

$$Z_i(x) = \frac{1}{h^d} K\left(\frac{X_i - x}{h}\right) - p(x).$$

The expectation of $Z_i(x_0)$ is the pointwise bias of a kernel density estimator with bandwidth $h$. Under the assumptions of the theorem, it is controlled in the following way:

$$|\mathbb{E} Z_i(x_0)| \leq L c_5 h^\beta.$$

Indeed,

$$\begin{aligned}
|\mathbb{E} Z_i(x_0)| &= \left|\int \frac{1}{h^d} K\left(\frac{t}{h}\right)\left[p(x_0 + t) - p(x_0)\right] dt\right| \\
&= \left|\int K(t)[p(x_0 + ht) - p(x_0)]\, dt\right| \\
&= \left|\int K(t)[p(x_0 + ht) - p_{x_0}^{(\beta)}(x_0 + ht)]\, dt \right. \\
&\quad \left. + \int K(t)[p_{x_0}^{(\beta)}(x_0 + ht) - p(x_0)]\, dt\right|.
\end{aligned} \qquad (A.6)$$

To control the first term in the right-hand side of (A.6), note that since $K$ has support $[-1,1]^d$, for any $h < r/\sqrt{d}$, we have $x_0 + ht \in \mathcal{B}(x_0, r)$ for any $t \in [-1,1]^d$. Thus, using the fact that $p$ is in $\Sigma(\beta, L, r, x_0)$, we have

$$\left|\int K(t)[p(x_0 + ht) - p_{x_0}^{(\beta)}(x_0 + ht)]\, dt\right| \leq L \int |K(t)| \|ht\|^\beta\, dt.$$

Now, since $K$ is a $\lfloor\beta\rfloor$-valid kernel (cf. Proposition A.2) and $p_{x_0}^{(\beta)} - p(x_0)$ is a polynomial of degree at most $\lfloor\beta\rfloor$ with no constant term, the second term in the right-hand side of (A.6) is zero. Therefore, we have

$$|\mathbb{E}Z_i(x_0)| \leq Lh^\beta \int |K(t)| \|t\|^\beta \, dt \qquad \text{for any } h \leq r.$$

Now denote, for simplicity, $Z_i = Z_i(x_0)$ and let $\overline{Z_i}$ be the centered version of $Z_i$. When $Lc_5 h^\beta \leq \delta/2$, we then have

$$\mathbb{P}\{|\hat{p}_n(x_0) - p(x_0)| \geq \delta\} \leq \mathbb{P}\left\{\frac{1}{n}\left|\sum_{i=1}^n \overline{Z_i}\right| \geq \delta - Lc_5 h^\beta\right\}$$

$$\leq \mathbb{P}\left\{\frac{1}{n}\left|\sum_{i=1}^n \overline{Z_i}\right| \geq \frac{\delta}{2}\right\}.$$

The right-hand side of the last inequality can be bounded by applying Bernstein's inequality (see Devroye *et al.* (1996), Theorem 8.4, page 124) to $\overline{Z_i}$ and $-\overline{Z_i}$ successively. For $h \leq 1$, we have

$$|\overline{Z_i}| \leq \|K\|_\infty h^{-d} + L^* + Lc_5 h^\beta \leq c_7 h^{-d},$$

where $c_7 = \|K\|_\infty + L^* + Lc_5$ and

$$\mathrm{Var}\{Z_i\} \leq h^{-d} \int K(u)^2 p(x_0 + hu) \, du \leq c_8 h^{-d},$$

where $c_8 = L^* \|K\|^2$. Applying Bernstein's inequality now yields

$$\mathbb{P}\{|\hat{p}_n(x_0) - p(x_0)| \geq \delta\} \leq 2\exp\left(-\frac{n(\delta/2)^2}{2(c_8 h^{-d} + c_7 h^{-d}\delta/6)}\right)$$

$$\leq 2\exp(-c_6 n h^d \delta^2)$$

for any $\delta \leq \Delta$ and where $\Delta = 6c_8/c_7$ and $c_6 = 1/(16c_8)$. □

### A.3. Equivalent formulation for the $\gamma$-exponent condition

The following proposition gives an equivalent formulation for the $\gamma$-exponent condition.

**Proposition A.1.** *Fix $\lambda > 0, \gamma > 0$ and $L_Q > 0$.*
  *Define $\mathcal{L} = \mathcal{L}(\lambda) = \{p = \lambda\}$. The two following statements are equivalent:*

(i) $\exists c > 0$ *and $\varepsilon_0 > 0$ such that for any $0 < \varepsilon \leq \varepsilon_0$, we have*

$$Q\{x \in \mathcal{X} : 0 < |p(x) - \lambda| \leq \varepsilon\} \leq c\varepsilon^\gamma;$$

(ii) $\exists c' > 0$ and $\varepsilon_1 > 0$ such that for any $0 < \varepsilon \leq \varepsilon_1$, we have

$$Q\{x \in \mathcal{X} : 0 < |p(x) - \lambda| \leq \varepsilon\} \leq L_Q,$$

and for all $G \subseteq \mathcal{X} \setminus \mathcal{L}$ satisfying $Q(G) \leq L_Q$, we have

$$Q(G) \leq c' \left( \int_G |p(x) - \lambda| \, dQ(x) \right)^{\gamma/(1+\gamma)}. \tag{A.7}$$

**Proof.** The proof of (i) $\Rightarrow$ (ii) essentially follows that of Tsybakov (2004), Proposition 1. Define

$$\varepsilon_1 = \varepsilon_0 \wedge \left( \frac{L_Q}{c(1+\gamma)} \right)^{1/\gamma}.$$

Observe that for any $0 < \varepsilon \leq \varepsilon_1$, we have

$$Q\{x \in \mathcal{X} : 0 < |p(x) - \lambda| \leq \varepsilon\} \leq c\varepsilon^\gamma \leq c\varepsilon_1^\gamma = \frac{L_Q}{1+\gamma} \leq L_Q.$$

Define $\mathcal{A}_\varepsilon = \{x : |p(x) - \lambda| > \varepsilon\}$ for all $0 < \varepsilon \leq \varepsilon_0$. For any measurable set $G \subset \mathcal{X} \setminus \mathcal{L}$, we have

$$\int_G |p(x) - \lambda| \, dQ(x) \geq \varepsilon Q(G \cap \mathcal{A}_\varepsilon)$$

$$\geq \varepsilon [Q(G) - Q(\mathcal{A}_\varepsilon^c \cap \mathcal{L}^c)]$$

$$\geq \varepsilon [Q(G) - \underline{c}\varepsilon^\gamma] \qquad \forall \underline{c} > c,$$

where the last inequality is obtained using (i). Maximizing the last term with respect to $\varepsilon > 0$, we get

$$\left( \int_G |p(x) - \lambda| \, dQ(x) \right)^{\gamma/(1+\gamma)} \geq Q(G) \left( \frac{\gamma}{1+\gamma} \right)^{\gamma/(1+\gamma)} \left( \frac{1}{1+\gamma} \right)^{1/(1+\gamma)} \underline{c}^{-1/(1+\gamma)}.$$

This yields (A.7) with $c' = e^{-2/e} \underline{c}^{1/(1+\gamma)}$. Note that the maximum is obtained for $\varepsilon = (\frac{Q(G)}{\underline{c}(1+\gamma)})^{1/\gamma} \leq \varepsilon_0$ for sufficiently large $\underline{c}$ and (i) is valid for this particular $\varepsilon$.

We now prove that (ii) $\Rightarrow$ (i). Consider $\varepsilon_1 > 0$ such that $Q(\mathcal{A}_\varepsilon^c \cap \mathcal{L}^c) \leq L_Q$ for any $0 < \varepsilon \leq \varepsilon_1$ and $c' > 0$ such that (A.7) is satisfied for any $G \subseteq \mathcal{X} \setminus \mathcal{L}$, $Q(G) \leq L_Q$. Taking $G = \mathcal{A}_\varepsilon^c \cap \mathcal{L}^c$ in (A.7) yields

$$Q\{x : 0 < |p(x) - \lambda| \leq \varepsilon\} = Q(\mathcal{A}_\varepsilon^c \cap \mathcal{L}^c)$$

$$\leq c' \left( \int_{\mathcal{A}_\varepsilon^c \cap \mathcal{L}^c} |p(x) - \lambda| \, dQ(x) \right)^{\gamma/(1+\gamma)}$$

$$\leq c' (\varepsilon Q(\mathcal{A}_\varepsilon^c \cap \mathcal{L}^c))^{\gamma/(1+\gamma)}.$$

Therefore,
$$Q\{x : 0 < |p(x) - \lambda| \leq \varepsilon\} \leq (c')^{1+\gamma}\varepsilon^{\gamma}.$$

This inequality yields (i) with $\varepsilon_0 = \varepsilon_1$ and $c = (c')^{1+\gamma}$. □

## A.4. On $\beta$-valid kernels

We recall here the definition of $\beta$-valid kernels and state a property that is useful in the present study.

**Definition A.1.** *Let $K$ be a real-valued function on $\mathbb{R}^d$, with support $[-1,1]^d$. For fixed $\beta > 0$, the function $K(\cdot)$ is said to be a $\beta$-valid kernel if it satisfies $\int K = 1$, $\int |K|^p < \infty$ for any $p \geq 1$, $\int \|t\|^{\beta}|K(t)|\,\mathrm{d}t < \infty$ and, in the case $\lfloor \beta \rfloor \geq 1$, it satisfies $\int t^s K(t)\,\mathrm{d}t = 0$ for any $s = (s_1, \ldots, s_d) \in \mathbb{N}^d$ such that $1 \leq s_1 + \cdots + s_d \leq \lfloor \beta \rfloor$.*

**Example A.1.** Let $\beta > 0$. For any $\beta$-valid kernel $K$ defined on $\mathbb{R}^d$, consider the product kernel
$$\tilde{K}(x) = K(x_1)K(x_2)\cdots K(x_d)\mathbb{1}_{x \in [-1,1]^d}$$
for any $x = (x_1, \ldots, x_d) \in \mathbb{R}^d$. It can then be easily shown that $\tilde{K}$ is a $\beta$-valid kernel on $\mathbb{R}^d$. Now, for any $\beta > 0$, an example of a 1-dimensional $\beta$-valid kernel is given in Tsybakov (2009), Section 1.2.2, the construction of which is based on Legendre polynomials. This eventually proves the existence of a multivariate $\beta$-valid kernel for any given $\beta > 0$.

The following proposition holds.

**Proposition A.2.** *Fix $\beta > 0$. If $K$ is a $\beta$-valid kernel, then $K$ is also a $\beta'$-valid kernel for any $0 < \beta' \leq \beta$.*

**Proof.** Fix $\beta$ and $\beta'$ such that $0 < \beta' \leq \beta$. Observe that $\lfloor \beta' \rfloor \leq \lfloor \beta \rfloor$ yields that if $\lfloor \beta' \rfloor \geq 1$, then for any $\beta$-valid kernel $K$, we have $\int t^s K(t)\,\mathrm{d}t = 0$ for any $s = (s_1, \ldots, s_d)$ such that $1 \leq s_1 + \cdots + s_d \leq \lfloor \beta' \rfloor$. It remains to check that

$$\int_{\mathbb{R}^d} \|t\|^{\beta'}|K(t)|\,\mathrm{d}t < \infty. \tag{A.8}$$

Consider the decomposition
$$\int_{\mathbb{R}^d} \|t\|^{\beta'}|K(t)|\,\mathrm{d}t = \int_{\|t\|\leq 1} \|t\|^{\beta'}|K(t)|\,\mathrm{d}t + \int_{\|t\|\geq 1} \|t\|^{\beta'}|K(t)|\,\mathrm{d}t$$
$$\leq \int_{\mathbb{R}^d} |K(t)|\,\mathrm{d}t + \int_{\|t\|\geq 1} \|t\|^{\beta}|K(t)|\,\mathrm{d}t.$$

To prove (A.8), note that since $K$ is a $\beta$-valid kernel, we have $\int_{\mathbb{R}^d} |K(t)|\,\mathrm{d}t < \infty$ and

$$\int_{\|t\|\geq 1} \|t\|^\beta |K(t)|\,\mathrm{d}t \leq \int_{\mathbb{R}^d} \|t\|^\beta |K(t)|\,\mathrm{d}t < \infty. \qquad \square$$

### A.5. Technical lemmas for minimax lower bounds

We collect here technical results that are used in Section 5. For a recent survey on the construction of minimax lower bounds, see Tsybakov (2009), Chapter 2. We first give a lemma related to subset extraction.

Fix an integer $k \geq 1$, and for any $\omega = (\omega_1, \ldots, \omega_k)$ and $\omega' = (\omega'_1, \ldots, \omega'_k)$ in $\{-1,1\}^k$ or in $\{0,1\}^k$, define the *Hamming distance* between $\omega$ and $\omega'$ by

$$\rho(\omega, \omega') = \sum_{i=1}^k \mathbb{1}_{\{\omega_i \neq \omega'_i\}}.$$

The following lemma can be found in Rigollet (2006), Lemma A.2. It is a straightforward corollary of Birgé and Massart (2001), Lemma 4, stated in a way which is more adapted to our purposes.

For any integers $N, \ell \geq 1$, define

$$\Omega_\ell = \left\{\omega \in \{0,1\}^N, \sum_{j=1}^N \omega_j = \ell\right\}$$

or, equivalently,

$$\Omega_\ell = \left\{\omega \in \{-1,1\}^N, \sum_{j=1}^N (\omega_j + 1)/2 = \ell\right\}.$$

**Lemma A.1 (Birgé and Massart (2001)).** *Let $N$ and $\ell$ be two integers such that $N \geq 6\ell \geq 6$. There then exists a subset $\Omega$ of $\Omega_{2\ell}$ such that*

$$\sum_{j=1}^N \mathbb{1}(\omega_j \neq \omega'_j) \geq \ell + 1 \qquad \forall \omega \neq \omega' \in \Omega,$$

*and $s = \mathrm{card}(\Omega)$ satisfies*

$$\log(s) \geq C\ell \log(N/\ell)$$

*for some numerical constant $C > 0$.*

The next lemma can be found in Tsybakov (2009), Theorem 2.7, and is stated here in a form adapted to the DLSE framework. It involves, in particular, the Kullback–Leibler divergence between two probability densities $p$ and $q$ on $\mathbb{R}^d$, defined by

$$K(p,q) = \begin{cases} \int_{\mathbb{R}^d} \log\left(\frac{p(x)}{q(x)}\right) p(x)\, \mathrm{d}x, & \text{if } P_p \ll P_q, \\ +\infty, & \text{otherwise.} \end{cases}$$

**Lemma A.2.** *Let $d$ be a pseudo-distance between subsets of $\mathcal{X} \subset \mathbb{R}^d$. Let $\mathcal{P}$ be a set of densities and assume that there exists a finite subset $\mathcal{N} \subset \mathcal{P}$ with $2 \leq \mathrm{card}(\mathcal{N}) = s < \infty$ such that*

$$d(\Gamma_p(\lambda), \Gamma_q(\lambda)) \geq 2\varepsilon \qquad \forall p, q \in \mathcal{N}, p \neq q, \tag{A.9}$$

*and*

$$\max_{q \in \mathcal{N}} K(q, p_0) \leq C \log(s) \tag{A.10}$$

*for some $p_0 \in \mathcal{N}$.*

*There then exists an absolute positive constant $C'$ such that for any estimator $\hat{G}_n$ of $\Gamma_p(\lambda)$ constructed from the sample $X_1, \ldots, X_n$, we have*

$$\sup_{p \in \mathcal{P}} \mathbb{E}[d(\Gamma_p(\lambda), \hat{G}_n)] \geq C'\varepsilon.$$

## Acknowledgements

The authors are thankful to the anonymous referees for their precious remarks that helped to improve the clarity of the text.